\documentclass{sig-alternate}
\usepackage{amsfonts,latexsym,amssymb,amsmath,algorithmic,epsfig,rotating}
\newlength{\jmr}
\newlength{\bjorn}
\renewcommand{\mod}{\mathrm{mod}}
\newtheorem{agpthm}{AGP Theorem}

\newcommand{\sat}{\mathtt{3CNFSAT}}

\newtheorem{dfn}{Definition}[section]
\newtheorem{algor}[dfn]{Algorithm}
 
\newtheorem{rem}[dfn]{Remark}

\newtheorem{prop}[dfn]{Proposition}
\newtheorem{thm}[dfn]{Theorem} 
\newtheorem{lemma}[dfn]{Lemma}
\newtheorem{cor}[dfn]{Corollary}
\newtheorem{ex}[dfn]{Example}

\newcommand{\thth}{^{\text{\underline{th}}}}
\newcommand{\rd}{^{\text{\underline{rd}}}}

\newcommand{\nd}{^{\text{\underline{nd}}}}

\newcommand{\ord}{{\mathrm{ord}}}

\newcommand{\np}{{\mathbf{NP}}}

\newcommand{\feas}{{\text{{\tt FEAS}}}}

\newcommand{\ffq}{{\feas_{\F_{\stackrel{\text{\scalebox{.8}[.8]{prime}}} 
{\text{\scalebox{.8}[.8]{powers}}}}}}}  
\newcommand{\fqp}{{\feas_{\Q_\mathrm{primes}}}}

\newcommand{\zpp}{{\mathbf{ZPP}}}

\newcommand{\pp}{\mathbf{P}}

\newcommand{\expt}{\mathbf{EXPTIME}}

\newcommand{\eps}{\varepsilon}
\newcommand{\Pro}{{\mathbb{P}}}

\newcommand{\F}{\mathbb{F}}
\newcommand{\Q}{\mathbb{Q}}
\newcommand{\R}{\mathbb{R}}
\newcommand{\C}{\mathbb{C}}
\newcommand{\N}{\mathbb{N}}
\newcommand{\Z}{\mathbb{Z}}

\newcommand{\cA}{{\mathcal{A}}}

\newcommand{\cE}{{\mathcal{E}}}

\newcommand{\cD}{{\mathcal{D}}}
\newcommand{\cF}{{\mathcal{F}}}

\newcommand{\cI}{{\mathcal{I}}}

\newcommand{\cL}{{\mathcal{L}}}

\newcommand{\cP}{{\mathcal{P}}}
\newcommand{\cR}{{\mathcal{R}}}
\newcommand{\cS}{{\mathcal{S}}}

\renewcommand{\qed}{$\blacksquare$}
\newcommand{\dia}{$\diamond$}

\newcommand{\newt}{\mathrm{Newt}}

\newcommand{\size}{\mathrm{size}}

\newcommand{\res}{\mathrm{Res}}
\newcommand{\supp}{\mathrm{Supp}}

\newlength{\hwl}
\begin{document}

\conferenceinfo{ISSAC'10,} {July 25--28, 2010, Munich, Germany.} 
\CopyrightYear{2010}
\crdata{xxx-x-xxxxx-xxx-x/xx/xx}

\title{Near $\text{\scalebox{2}[2]{$\np$}}$-Completeness for Detecting 
$\text{\scalebox{2}[2]{$p$}}$-adic Rational Roots in One Variable} 
\numberofauthors{3} 
\author{ 
\alignauthor Martin Avenda\~{n}o$^{^{^\text{{\normalsize $*$}}}}$\\
\affaddr{TAMU 3368}\\ 
\affaddr{Mathematics Dept.}\\
\affaddr{\scalebox{.6}[1]{College Station, TX \ 77843-3368, USA}}
\email{\scalebox{.7}[1]{mavendar@yahoo.com.ar}} 
\alignauthor Ashraf Ibrahim$^{^{^\text{{\normalsize $*$}}}}$\\
\affaddr{TAMU 3368}\\ 
\affaddr{Mathematics Dept.}\\
\affaddr{\scalebox{.6}[1]{College Station, TX \ 77843-3368, USA}}
\email{\scalebox{.7}[1]{aibrahim@math.tamu.edu}} 
\alignauthor
J. Maurice Rojas\titlenote{
Partially supported by NSF individual grant DMS-0915245 and NSF CAREER 
grant DMS-0349309. Rojas was also partially supported by Sandia National 
Laboratories.}\\ 
\affaddr{TAMU 3368}\\ 
\affaddr{Mathematics Dept.}\\
\affaddr{\scalebox{.6}[1]{College Station, TX \ 77843-3368, USA}}
\email{\scalebox{.7}[1]{rojas@math.tamu.edu}} 
\alignauthor
Korben Rusek$^{^{^\text{{\normalsize $*$}}}}$\\
\affaddr{TAMU 3368}\\ 
\affaddr{Mathematics Dept.}\\
\affaddr{\scalebox{.6}[1]{College Station, TX \ 77843-3368, USA}}
\email{\scalebox{.7}[1]{korben@rusek.org}}}

\date{\today} 

\maketitle

\begin{abstract} 
We show that deciding whether a sparse univariate polynomial  
has a $p$-adic rational root can be done in $\np$ for 
most inputs. We also prove a polynomial-time upper bound 
for trinomials with suitably generic $p$-adic Newton polygon. 
We thus improve the best previous complexity upper bound of 
$\expt$. We also prove an unconditional complexity 
lower bound of $\np$-hardness with respect to randomized reductions 
for general univariate polynomials. The best previous lower bound assumed an 
unproved hypothesis 
on the distribution of primes in arithmetic progression.  We also discuss how 
our results complement analogous results over the real numbers. 
\end{abstract} 

\section{Introduction}  
The fields $\R$ and $\Q_p$ (the reals and the $p$-adic rationals) bear more
in common than just completeness with
respect to a metric: increasingly, complexity results for one field
have inspired and motivated analogous results in the other 
(see, e.g., \cite{cohenqe,vandenef} and the pair of 
works \cite{few} and \cite{amd}). We continue this theme by transposing recent 
algorithmic results for sparse polynomials over the real numbers \cite{brs} 
to the $p$-adic rationals, sharpening the underlying complexity bounds 
along the way (see Theorem \ref{thm:qp} below).  

More precisely, for any commutative ring $R$ with multiplicative identity, 
we let $\feas_R$ --- the {\bf $R$-feasibility}\linebreak 
\scalebox{.95}[1]{{\bf problem} (a.k.a.\ 
Hilbert's Tenth Problem over $R$ \cite{h10})}\linebreak 
--- denote the problem of deciding whether an input polynomial system 
$F\!\in\!\bigcup_{k,n\in\N} (\Z[x_1,\ldots,x_n])^k$ has a root in $R^n$. 
(The underlying input size is clarified in Definition \ref{dfn:basic} 
below.) Observe that $\feas_\R$, $\feas_\Q$, and 
$\{\feas_{\F_q}\}_{q \text{ a prime power}}$ are central 
problems respectively in algorithmic real algebraic geometry, algorithmic 
number theory, and cryptography. 

In particular, for any prime $p$ and $x\!\in\!\Z$, recall that the 
{\bf $p$-adic valuation}, $\ord_p x$, is the greatest $k$ such that $p^k|x$. 
We can extend $\ord_p(\cdot)$ to $\Q$ by
$\ord_p\left(\frac{a}{b}\right)\!:=\!\ord_p(a)-\ord_p(b)$ for any
$a,b\!\in\!\Z$; and we let $|x|_p\!:=\!p^{-\ord_p x}$ denote the {\bf
$p$-adic norm}. The norm $|\cdot|_p$ defines a natural metric satisfying
the ultrametric inequality and $\Q_p$ is, to put it tersely, the 
completion of $\Q$ with respect to this metric. This metric, 
along with $\ord_p(\cdot)$, extends naturally to the {\bf $p$-adic
complex numbers} $\C_p$, which is the metric completion of the algebraic 
closure of $\Q_p$ \cite[Ch.\ 3]{robert}.

We will also need to recall the following containments of 
complexity classes: 
$\pp\!\subseteq\!\zpp\!\subseteq\!\np\!\subseteq\cdots\subseteq\!\expt$, 
and the fact that the properness of {\bf every}  
inclusion above (save $\pp\!\subsetneqq\!\expt$) is a major open problem 
\cite{lab,papa}. The definitions of the aforementioned 
complexity classes are reviewed briefly in the Appendix (see 
also \cite{papa} for an excellent textbook treatment). 

\subsection{The Ultrametric Side: Relevance and\\ Results} 
\label{sub:padic} 
Algorithmic results over the $p$-adics are central in 
many computational areas: polynomial time factoring algorithms over $\Q[x_1]$ 
\cite{lll}, computational complexity \cite{antsv}, studying prime ideals in 
number fields \cite[Ch.\ 4 \& 6]{cohenant}, elliptic 
curve cryptography \cite{lauder}, and the computation of zeta functions 
\cite{denefver}. Also, much work has gone into using
$p$-adic methods to algorithmically detect
rational points on algebraic plane curves via variations of the {\bf Hasse
Principle}\footnote{
If $F(x_1,\ldots,x_n)\!=\!0$ is any polynomial equation and
$Z_K$ is its zero set in $K^n$, then the Hasse Principle is the
assumption that [$Z_\C$ smooth, $Z_\R\!\neq\!\emptyset$, and
$Z_{\Q_p}\!\neq\!\emptyset$ for all primes $p$] implies
$Z_\Q\!\neq\!\emptyset$ as well. The Hasse Principle is a theorem when
$Z_\C$ is a quadric hypersurface or a curve of genus zero,
but fails in subtle ways already for curves of genus one (see, e.g.,
\cite{bjornhasse1}). } (see, e.g., \cite{colliot,bjornhasse2,bjornbm}). 
However, our knowledge of the complexity of deciding the existence of 
solutions for {\bf sparse} polynomial equations over $\Q_p$ is surprisingly 
coarse: good bounds for the number of solutions over $\Q_p$ in 
one variable weren't even known until the late 1990s \cite{lenstra2}. 
So we focus on precise complexity bounds for one variable. 
\begin{dfn}
\label{dfn:basic}
Let $f(x)\!:=\!\sum^m_{i=1} c_i
x^{a_i}\!\in\!\Z[x_1,\ldots,x_n]$\linebreak
where $x^{a_i}\!:=\!x^{a_{1i}}_1\cdots x^{a_{ni}}_n$, $c_i\!\neq\!0$ for all
$i$, and the $a_i$ are pair-wise distinct. We call such an $f$ an
{\bf $\pmb{n}$-variate $\pmb{m}$-nomial}.
Let us also define\\
\mbox{}\hfill $\size(f)\!:=\!\sum^m_{i=1}
\log_2\left[(2+|c_i|)(2+|a_{1,i}|)\cdots
(2+|a_{n,i}|)\right]$\hfill\mbox{}\\
and, for any
$F\!:=\!(f_1,\ldots,f_k)\!\in\!(\Z[x_1,\ldots,x_n])^k$,
we\linebreak
define $\size(F)\!:=\!\sum^k_{i=1}\size(f_i)$.
Finally, we let $\cF_{n,m}$ denote the subset of $\Z[x_1,\ldots,x_n]$
consisting of polynomials with\linebreak exactly $m$ monomial terms
\dia
\end{dfn}
For instance, $\size(1+cx^{99}_1+x^d_1)\!=\!\Theta(\log(c)+\log(d))$. So
the degree, $\deg f$, of a polynomial $f$ can sometimes be exponential
in its size. Note also that $\Z[x_1]$ is the disjoint
union $\bigsqcup_{m\geq 0}\cF_{1,m}$.
\begin{dfn}
\label{dfn:qp}
Let $\fqp$ denote the problem of deciding, for an input polynomial
system $F$\linebreak
$\in\!\bigcup_{k,n\in\N}(\Z[x_1,\ldots,x_n])^k$ {\bf and}
an input prime $p$, whether $F$ has a root in $\Q^n_p$. 
Also let $\Pro\!\subset\!\N$ denote the set of primes and, when 
$\cI$ is a family of such pairs $(F,p)$, we let $\fqp(\cI)$ denote the 
restriction of $\fqp$ to inputs in $\cI$. The underlying input sizes 
for $\fqp$ and $\fqp(\cI)$ shall be $\size_p(F)\!:=\!\size(F)+\log p$ (cf.\ 
Definition \ref{dfn:basic}). Finally, let $(\Z\times (\N\cup \{0\}))^\infty$
denote the set of all infinite sequences of pairs $((c_i,a_i))^\infty_{i=1}$ 
with $c_i\!=\!a_i\!=\!0$ for $i$ sufficiently large. \dia
\end{dfn}
\begin{rem}
Note that $\Z[x_1]$ admits a natural embedding
into $(\Z\times (\N\cup \{0\}))^\infty$ by considering coefficient-exponent
pairs in order of increasing exponents, e.g., $a+bx^{99}+x^{2001} \mapsto
((a,0),(b,99),(1,2001),(0,0),(0,0),\ldots)$. \dia
\end{rem}

While there are 
now randomized algorithms for factoring $f\!\in\!\Z[x_1]$ over $\Q_p[x_1]$ with 
expected complexity polynomial in $\size_p(f)+\deg(f)$ 
\cite{cantorqp} (see also \cite{chistov}), no such algorithms are known   
to have complexity polynomial in $\size_p(f)$ alone. Our main theorem below 
shows that such algorithms are hard to find 
because their existence is essentially equivalent to the 
$\pp\!=\!\np$ problem. Moreover, we obtain new sub-cases
of $\fqp(\Z[x_1]\times\Pro)$ lying in $\pp$. 
\begin{thm}
\label{thm:qp} \mbox{}\\
1.\ $\fqp(\cF_{1,k}\times \Pro)\!\in\!\pp$ for $k\!\in\!\{0,1,2\}$.\\ 
2.\ For any $f(x_1)\!=\!c_1+c_2x^{a_2}_1+c_3x^{a_3}_1\!\in\!\Z[x_1]$ with 
the points\linebreak
\mbox{}\hspace{.4cm}$\{(0,\ord_p(c_1)),(a_2,\ord_p(c_2)),(a_3,\ord_p(c_3))\}$ 
{\bf non}-collinear,\linebreak 
\mbox{}\hspace{.4cm}and $p$ {\bf not} dividing $a_2$, $a_3$, or $a_3-a_2$, 
we can decide the\linebreak
\mbox{}\hspace{.4cm}existence of a root in $\Q_p$ for $f$ in $\pp$.\\
3.\ There is a countable union of algebraic hypersurfaces\linebreak 
\mbox{}\hspace{.4cm}$E\!\subsetneqq\!\Z[x_1]\times \Pro$, with natural 
density $0$, such that\linebreak 
\mbox{}\hspace{.4cm}$\fqp((\Z[x_1]\times \Pro)\setminus E)\!\in\!\np$. 
Furthermore, we can\linebreak 
\mbox{}\hspace{.4cm}decide in $\pp$ whether an $f\!\in\!\cF_{1,3}$ 
also lies in $E$.\\
4.\ If $\fqp(\Z[x_1]\times \Pro)\!\in\!\zpp$ then $\np\!\subseteq\!\zpp$.\\
5.\ If the Wagstaff Conjecture is true, then 
$\fqp(\Z[x_1])$\linebreak
\mbox{}\hspace{.4cm}$\in\!\pp \Longrightarrow \pp\!=\!\np$, i.e., we can 
strengthen Assertion (4)\linebreak 
\mbox{}\hspace{.4cm}above. 
\end{thm} 
\begin{rem}
The Wagstaff Conjecture, dating back to 1979 
(see, e.g., \cite[Conj.\ 8.5.10, pg.\ 224]{bs}), is the\linebreak 
assertion that the least prime congruent to $k$ mod $N$ is\linebreak  
$O(\varphi(N)\log^2 N)$, 
where $\varphi(N)$ is the number of integers in $\{1,\ldots,N\}$ relatively 
prime to $N$. Such a bound is significantly stronger than the known 
implications of the {\bf Generalized Riemann Hypothesis (GRH)}. \dia 
\end{rem} 

While the real analogue of Assertion (1) is known (and easy), the 
stronger real analogue $\feas_\R(\cF_{1,3})\!\in\!\pp$ to 
Assertion (2) was unknown until \cite[Thm.\ 1.3]{brs}. We hope 
to strengthen Assertion (2) to $\fqp(\cF_{1,3}\times\Pro)\!\in\!\pp$ in 
future work. In fact, we can attain polynomial complexity already for more 
inputs in $\cF_{1,3}\times \Pro$ than stated above, and this is clarified in 
Section \ref{sec:qp}. 

Note that $\Q_p$ is uncountable and thus, unlike $\feas_{\F_p}$,
$\feas_{\Q_p}$ does {\bf not} admit an obvious succinct certificate. Indeed,
while it has been known since the late 1990's that $\fqp\!\in\!\expt$
relative to our notion of input size \cite{mw1,mw2},
we are unaware of any earlier algorithms
yielding $\fqp(\Z[x_1,\ldots,x_n]\times\Pro)\!\in\!\np$ for any fixed $n$:
even $\fqp(\cF_{1,4}\times\Pro)\text{\scalebox{1}[.85]{$\stackrel{?}{\in}$}}
\np$ and $\feas_\R(\cF_{1,4})\text{\scalebox{1}[.85]{$\stackrel{?}{\in}$}}\np$
are open questions.\footnote{An earlier result claiming
$\fqp(\Z[x_1]\times \Pro)\!\in\!\np$ for ``most'' inputs 
\cite[Main Thm.]{myqua} appears to have fatal errors in its proof.}
Practically speaking, zero density means that under most reasonable
input restrictions, the algorithmic speed-up in Assertion (3) is valid over a
significantly large fraction of inputs.
\begin{ex}
Let $T$ denote the family of pairs
$(f,p)\!\in\!\Z[x_1]\times \Pro$ with $f(x_1)\!=\!a+bx^{11}_1+cx^{17}_1
+x^{31}_1$ and let $T^*\!:=\!T\setminus E$. Then there is
a sparse $61\times 61$ structured matrix $\cS$ (cf.\ Lemma \ref{lemma:syl}
in Section \ref{sub:transfer} below), whose entries lie in
$\{0,1,31,a,b,11b,c,17c\}$,
such that $(f,p)\!\in\!T^* \Longleftrightarrow p\!\not|\!\det \cS$.
So by Theorem \ref{thm:qp}, $\fqp(T^*)\!\in\!\np$, and Corollary \ref{cor:lots}
in Section \ref{sec:qp} below tells us that for large coefficients,
$T^*$ occupies almost all of $T$. In particular,
letting $T(H)$ (resp.\ $T^*(H)$) denote those
pairs $(f,p)$ in $T$ (resp.\ $T^*$) with\linebreak 
\scalebox{.94}[1]{$|a|,|b|,|c|,p\!\leq\!H$, we have
$\frac{\#T^*(H)}{\#T(H)}\!\geq\!\left(1-\frac{61}{H}\right)
\left(1-\frac{31\log_2(124H)}{H}\right)$.}\linebreak  
For instance, one can check via {\tt Maple} that\\ 
\mbox{}\hfill 
$(-973+21x^{11}_1-2x^{17}_1 +x^{31}_1,p)\!\in\!T^*$\hfill\mbox{}\\ 
for all but $352$ primes $p$. \dia 
\end{ex}

\noindent 
The exceptions in Assertion (3) appear to be due to the presence of
{\bf ill-conditioned} polynomials: $f$ having a root $\zeta$ with
the ($p$-adic) norm of $f'(\zeta)$ very small --- a phenomenon of
approximation present in complete fields like $\R$, $\C$, and $\Q_p$. 
Curiously, the real analogue of Assertion (3) remains unknown 
\cite[Sec.\ 1.2]{brs}.  

As for lower bounds, while it is not hard to show that the 
full problem $\fqp$ is $\np$-hard from scratch, the least $n$ making 
$\fqp(\Z[x_1,\ldots,x_n]\times \Pro)$ $\np$-hard appears not to have been known 
unconditionally. In particular, a weaker version of Assertion (4) 
was found recently, but only under the truth of an unproved 
hypothesis on the distribution of primes in arithmetic progresion 
\cite[Main Thm.]{myqua}. Assertion (4) thus also provides an interesting 
contrast to earlier work of H.\ W.\ Lenstra, Jr.\ \cite{lenstra1}, who showed
that one can actually find all {\bf low} degree factors of a sparse polynomial 
(over $\Q[x_1]$ as opposed to $\Q_p[x_1]$) in polynomial time.

\subsection{Random Primes and Tropical Tricks} 
\label{sub:key} 
The key to proving our lower bound results (Assertions (4) and (5) of 
Theorem \ref{thm:qp}) is an efficient reduction from a 
problem discovered to be $\np$-hard by David Alan Plaisted: deciding whether  
a sparse univariate polynomial vanishes at a complex $D\thth$ root of unity 
\cite{plaisted,mega}. Reducing from this problem to its analogue over 
$\Q_p$ is straightforward, provided $\Q^*_p$ contains a cyclic subgroup 
of order $D$ where $D$ has sufficiently many distinct prime divisors. We thus  
need to consider the factorization of $p-1$, which in turn leads us to 
primes congruent to $1$ modulo certain integers. 

While efficiently constructing random primes in {\bf arbitrary} 
arithmetic progressions remains a famous open problem, we can now at least 
efficiently build random primes $p$ such that $p$ is moderately sized but 
$p-1$ has many prime factors. 
We use the notation $[j]\!:=\!\{1,\ldots,j\}$ for any $j\!\in\!\N$. 
\begin{thm}
\label{thm:von}
For any $\delta\!>\!0$, a failure probability\linebreak 
$\eps\!\in\!(0,1/2)$, and $n\!\in\!\N$, we can find --- 
within\linebreak 
$O\!\left((n/\eps)^{\frac{3}{2}+\delta} + 
\left(n\log(n)+\log\frac{1}{\eps}\right)^{7+\delta}
\right)$ randomized bit\linebreak 
operations --- a sequence $P\!=\!(p_i)^n_{i=1}$ of consecutive primes 
and a positive integer $c$ such that\\ 
\mbox{}\hfill $\log(c),\log\left(\prod\limits^n_{i=1}p_i\right) = 
O(n\log(n)+\log(s/\eps))$\hfill\mbox{}\\ 
and, with probability $\geq\!1-\eps$, the number 
$p\!:=\!1+c\prod\limits^n_{i=1} p_i$ is prime. 
\end{thm}

\noindent
Theorem \ref{thm:von} and its proof are inspired in large part by an algorithm
of von zur Gathen, Karpinski, and Shparlinski \cite[Algorithm following Fact
4.9]{von}.  In particular, they used an intricate random sampling
technique \cite[Thm.\ 4.10]{von} to show, in our notation, that
the enumerative analogue of $\ffq(\Z[x_1,x_2])$ is $\#\pp$-hard
\cite[Thm.\ 4.11]{von}. Note in particular that neither of Theorem 4.10 of 
\cite{von} or Theorem \ref{thm:von} above implies the other. 

Our harder upper bound results (Assertions (2) and (3) of 
Theorem \ref{thm:qp}) will follow from 
an arithmetic analogue of toric deformations. Here, this simply means 
that we find ways to reduce problems involving general $f\!\in\!\Z[x_1]$ to 
similar problems involving binomials. As a warm-up, let us 
recall that the convex hull of any subset
$S\!\subseteq\!\R^2$ is the smallest convex set containing $S$. Also,
an edge of a polygon $P\!\subset\!\R^2$ is called {\bf lower} iff it has
an inner normal with positive last coordinate, and the {\bf lower hull} of
$P$ is simply the union of all its lower edges. 
\begin{lemma}
\label{lemma:newt}
(See, e.g., \cite[Ch.\ 6, sec.\ 1.6]{robert}.)
Given any polynomial $f(x_1)\!:=\!\sum^m_{i=1}c_ix^{a_i}_1\!\in\!\Z[x_1]$,
we define its {\bf $p$-adic Newton polygon},
$\newt_p(f)$, to be the convex hull of the points
$\{(a_i,\ord_p c_i)\; | \; i\!\in\!\{1,\ldots,m\}\}$. Then the number of roots
of $f$ in $\C_p$ with valuation $v$, counting multiplicities, is {\bf exactly} 
the horizontal length of the lower face of $\newt_p(f)$ with inner normal 
$(v,1)$. \qed
\end{lemma} 
\begin{ex}
For the polynomial\\ 
$f(x_1)\!:=\!243x^6-3646x^5+18240x^4-35310x^3+29305x^2
-8868x+36$,
the polygon $\newt_3(f)$ can easily be verified to resemble the
following illustration: \\
\mbox{}\hfill\epsfig{file=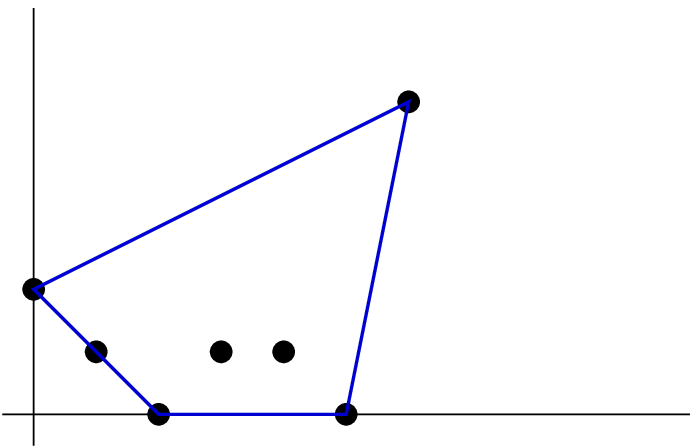,height=1.5in} \hfill\mbox{}\\
Note in particular that there are exactly $3$ lower edges, and
their respective horizontal lengths and inner normals are
$2$, $3$, $1$, and $(1,1)$, $(0,1)$, and $(-5,1)$. Lemma \ref{lemma:newt}
then tells us that $f$ has exactly $6$ roots in $\C_3$:
$2$ with $3$-adic valuation $1$, $3$ with $3$-adic valuation $0$, and $1$
with $3$-adic valuation $-5$. Indeed, one can check that the roots
of $f$ are exactly $6$, $1$, and $\frac{1}{243}$, with respective
multiplicities $2$, $3$, and $1$. \dia
\end{ex}

The binomial associated to summing the terms of $f$ corresponding to 
the vertices of a lower edge of $\newt_p(f)$ containing no other point of 
the form $(a_i,\ord_p c_i)$ in its interior is called a {\bf lower binomial}. 
\begin{lemma}
\label{lemma:ai} 
Suppose $f(x_1)\!=\!c_1+c_2x^{a_2}_1+c_3x^{a_3}_1\!\in\!\Z[x]$, 
the points $\{(0,\ord_p(c_1)),(a_2,\ord_p(c_2)),(a_3,\ord_p(c_3))\}$ 
are\linebreak  
{\bf non}-collinear, and $p$ is a prime {\bf not} dividing 
$a_2$, $a_3$, or $a_3-a_2$. Then the number of roots of 
$f$ in $\Q_p$ is exactly the number of roots of the 
$p$-adic lower binomials of $f$ in $\Q_p$. \qed  
\end{lemma} 

\noindent 
Our last lemma follows easily (taking direct limits) from a more general 
result (\cite[Thm.\ 4.5]{ai}) 
relating the number of roots of $f$ with the number of roots of 
its lower binomials over $\Z/p^N\Z$ for $N$ sufficiently large.

Our main results are proved in Section \ref{sec:qp}, after
the development of some additional theory below. 

\section{Background and Ancillary\\ Results} 
\label{sec:back}
Our lower bounds will follow from a common chain of reductions, 
so we will begin by reviewing the fundamental problem from which we reduce. 
We then show how to efficiently construct random primes 
$p$ such that $p-1$ has many prime factors 
in Section \ref{sub:agp},  
and conclude with some quantitative results for 
transferring complexity results over $\C$ to $\Q_p$ 
in Section \ref{sub:transfer}. 

\subsection{Roots of Unity and NP-Completeness}  
\label{sub:cyclo} 
Recall that any Boolean expression of one of the following forms:\\ 
$(\heartsuit)$ $y_i\vee y_j \vee y_k$, \
$\neg y_i\vee y_j \vee y_k$, \
$\neg y_i\vee \neg y_j \vee y_k$,  \
$\neg y_i\vee \neg y_j \vee \neg y_k$,\\ 
\mbox{}\hspace{.6cm}with $i,j,k\!\in\![3n]$,\\
is a $\sat$ {\bf clause}. 
Let us first refine slightly Plaisted's elegant reduction from $\sat$ to 
feasibility testing for univariate polynomial systems over the complex 
numbers \cite[Sec.\ 3, pp.\ 127--129]{plaisted}. 
\begin{dfn}  
\label{dfn:plai} 
Letting $P\!:=\!(p_{1},\ldots,p_{n})$ denote any\linebreak 
strictly increasing sequence of primes, let us inductively define a semigroup 
homomorphism $\cP_P$ --- 
the {\bf Plaisted morphism with respect to $P$} --- from certain Boolean 
expressions in the variables $y_1,\ldots,y_n$ to $\Z[x_1]$, as 
follows:\footnote{
Throughout this paper, for Boolean expressions, we will always identify $0$
with ``{\tt False}'' and  $1$ with ``{\tt True}''.}
(0) $D_P\!:=\!\prod^n_{i=1}p_{i}$, 
(1) $\cP_P(0)\!:=\!1$, (2) $\cP_P(y_i)\!:=\!x^{D_P/p_{i}}_1-1$, 
(3) $\cP_P(\neg B):=$
$(x^{D_P}_1-1)/\cP_P(B)$, for any 
Boolean expression $B$ for which $\cP_P(B)$ has already been defined, 
(4) $\cP_P(B_1\vee B_2)\!:=\!\mathrm{lcm}(\cP_P(B_1),\cP_P(B_2))$, 
for any Boolean expressions $B_1$ and $B_2$ for which $\cP_P(B_1)$ and 
$\cP_P(B_2)$ have already been defined. \dia 
\end{dfn} 

\begin{lemma} 
\label{lemma:plai} 
\cite[Sec.\ 3, pp.\ 127--129]{plaisted} 
Suppose $P\!=\!(p_i)^n_{k=1}$ is an increasing sequence of 
primes with $\log(p_{k})\!=\!O(k^\gamma)$ for some constant $\gamma$. Then, 
for all $n\!\in\!\N$ and any clause $C$ of the form 
$(\heartsuit)$, we have $\size(\cP_P(C))$ polynomial in $n$. 
In particular, 
$\cP_P$ can be evaluated at any such $C$ in time polynomial in $n$. 
Furthermore, if $K$ is any field possessing $D_P$ distinct ${D_P}\thth$ 
roots of unity, then a $\sat$ instance $B(y)\!:=C_1(y)\wedge \cdots \wedge 
C_k(y)$ has a satisfying assignment iff the univariate polynomial system 
$F_B\!:=\!(\cP_P(C_1), \ldots,\cP_P(C_k))$ has a root $\zeta\!\in\!K$ 
satisfying $\zeta^{D_P}-1$. \qed 
\end{lemma} 

\noindent
Plaisted actually proved the special case $K\!=\!\C$ of the above lemma, in 
slightly different language, in \cite{plaisted}. 
However, his proof extends verbatim to the more general 
family of fields detailed above. 

\subsection{Randomization to Avoid Riemann\\ Hypotheses}
\label{sub:agp} 
The result below allows us to prove Theorem \ref{thm:von} and further 
tailor Plaisted's clever reduction to our purposes.  
We let $\pi(x)$ the number of primes $\leq\!x$, and 
let $\pi(x;M,1)$ denote the number of primes $\leq\!x$ that are congruent 
to $1 \; \mod \; M$. 
\begin{agpthm} (very special case of \cite[Thm.\ 2.1, pg.\ 712]{carmichael})  
There exist $x_0\!>\!0$ and an $\ell\!\in\!\N$ such that for each  
$x\!\geq\!x_0$, there is a subset $\cE(x)\!\subset\!\N$ of finite cardinality 
$\ell$ with the following property: If $M\!\in\!\N$ satisfies 
$M\!\leq\!x^{2/5}$ and $a\not|M$ for all $a\!\in\!\cE(x)$ then 
$\pi(x;M,1)\!\geq\!\frac{\pi(x)}{2\varphi(M)}$. \qed  
\end{agpthm} 

\noindent 
For those familiar with \cite[Thm.\ 2.1, pg.\ 712]{carmichael}, 
the result above follows immediately upon specializing the 
parameters there as follows:\\ 
\mbox{}\hfill $(A,\eps,\delta,y,a)\!=\!(49/20,1/2,2/245, x,1)$ \hfill 
\mbox{}\\ 
(see also \cite[Fact 4.9]{von}). 

The AGP Theorem enables us to construct random primes from 
certain arithmetic progressions with high probability. An additional 
ingredient that will prove useful is the famous recent {\bf AKS algorithm} 
for deterministic polynomial-time primality checking \cite{aks}. 
Consider now the following algorithm. 
\begin{algor}\mbox{}\\ 
\label{algor:primes}{\bf Input:} A constant $\delta\!>\!0$, a failure 
probability $\eps\!\in\!(0,1/2)$, a positive integer $n$, and the 
constants $x_0$ and $\ell$ from the AGP Theorem.\\ 
{\bf Output:} An increasing sequence $P\!=\!(p_j)^n_{j=1}$ 
of primes such that $\log p\!=\!O(n\log(n)+\log(1/\eps))$ and, 
with probability $1-\eps$,  $p\!:=\!1+c\prod^n_{i=1} p_i$ 
is prime. In particular, the output always gives a true declaration 
as to the primality of $p$. 

\medskip 
\noindent
{\bf Description:}
\vspace{-.2cm}
\begin{enumerate}
\setcounter{enumi}{-1}
\item{Let $L\!:=\!\lceil 2/\eps\rceil\ell$ 
and compute the first $nL$ primes $p_1, \ldots,$\linebreak 
$p_{nL}$ in increasing order.}  
\vspace{-.2cm} 
\item{Define (but do not compute) $M_j\!:=\!\prod\limits^{jn}_{k=(j-1)n+1} p_k$ 
for any $j\!\in\!\N$. Then compute $M_L$, $M_i$ for a uniformly random 
$i\!\in\![L]$, 
and $x\!:=\!\max\left\{x_0,
17,
1+M^{5/2}_L
\right\}$. } 
\vspace{-.2cm} 
\item{\scalebox{.97}[1]{Compute $K\!:=\!\lfloor (x-1)/M_i\rfloor$ 
and $J\!:=\!\lceil 2\log(2/\eps)\log x\rceil$.}} 
\vspace{-.2cm} 
\item{Pick uniformly random $c\!\in\![K]$ until 
one either has $p\!:=\!1+cM_i$ prime, or one has 
$J$ such numbers that are each composite (using primality checks 
via the AKS algorithm along the way). } 
\vspace{-.2cm} 
\item{If a prime $p$ was found then output\\ 
\mbox{}\hfill ``{\tt $1+c\prod^{in}_{j=(i-1)n+1}p_j$ 
is a prime that works!}''\hfill\mbox{}\\ 
and stop. Otherwise, stop and output\\
``{\tt I have failed to find a suitable prime. 
Please forgive me.}'' \dia} 
\end{enumerate}
\end{algor}
\begin{rem} 
\label{rem:tran} 
In our algorithm above, it suffices to find integer approximations to the 
underlying logarithms and square-roots. In particular, we restrict to 
algorithms that can compute the $\log_2 \cL$ most significant bits of 
$\log \cL$, and the $\frac{1}{2}\log_2 \cL$ most significant bits of 
$\sqrt{\cL}$, using\\
\mbox{}\hfill $O((\log \cL)(\log \log \cL)\log \log \log \cL)$ 
\hfill\mbox{}\\
bit operations. Arithmetic-Geometric Mean Iteration and (suitably tailored) 
Newton Iteration are algorithms that respectively satisfy our 
requirements (see, e.g., \cite{dan} for a detailed description). \dia 
\end{rem} 

\noindent 
{\bf Proof of Theorem \ref{thm:von}:} 
It clearly suffices to prove that Algorithm \ref{algor:primes} 
is correct, has a success probability that is at least $1-\eps$, and 
works within\\ 
\mbox{}\hfill $O\!\left(\left(\frac{n}{\eps}\right)^{\frac{3}{2}+\delta}+
(n\log(n)+\log(1/\eps))^{7+\delta}\right)$\hfill\mbox{}\\ 
randomized bit operations, 
for any $\delta\!>\!0$. These assertions are proved directly below. \qed  

\smallskip 
\noindent
{\bf Proving Correctness and the Success Probability\linebreak 
Bound for Algorithm \ref{algor:primes}:} 
First observe that $M_1,\ldots,M_{L}$ are relatively prime. So at most
$\ell$ of the $M_i$ will be divisible by elements of $\cE(x)$.
Note also that $K\!\geq\!1$ and
$1+cM_i\!\leq\!1+KM_i\!\leq\!1+((x-1)/M_i)M_i\!=\!x$
for all $i\!\in\![L]$ and $c\!\in\![K]$.

Since $x\!\geq\!x_0$ and
$x^{2/5}\!\geq\!(x-1)^{2/5}\!\geq\!\left(M^{5/2}_i\right)^{2/5}\!=\!M_i$ for
all $i\!\in\![L]$, the AGP Theorem implies that with probability
$\geq 1-\frac{\eps}{2}$ (since $i\!\in\![\lceil 2/\eps\rceil \ell]$ is
uniformly random), the arithmetic
progression $\{1+M_i,\ldots,1+KM_i\}$ contains at
least $\frac{\pi(x)}{2\varphi(M_i)}\!\geq\!\frac{\pi(x)}{2M_i}$ primes.
In which case, the proportion of numbers in $\{1+M_i,\ldots,1+KM_i\}$ that
are prime is $\frac{\pi(x)}{2KM_i}\!>\!\frac{\pi(x)}{2+2KM_i}\!>\!
\frac{x/\log x}{2x}\!=\!\frac{1}{2\log x}$, since
$\pi(x)\!>\!x/\log x$ for all $x\!\geq\!17$ \cite[Thm.\ 8.8.1, pg.\ 233]{bs}.
So let us now assume that $i$ is fixed and $M_i$ is not divisible by
any element of $\cE(x)$.

Recalling the inequality $\left(1-\frac{1}{t}\right)^{ct}\!\leq\!e^{-c}$
(valid for all $c\!\geq\!0$ and $t\!\geq\!1$), we then see that
the AGP Theorem implies that the probability of {\bf not} finding a prime
of the form $p\!=\!1+cM_i$ after picking $J$ uniformly random
$c\!\in\![K]$ is
$\left(1-\frac{1}{2\log x}\right)^J \!\leq\!\left(1-\frac{1}{2\log x}
\right)^{2\log(2/\eps)\log x}\!\leq\!e^{-\log(2/\eps) }\!=\!
\frac{\eps}{2}$.

In summary, with probability $\geq\!1-\frac{\eps}{2}
-\frac{\eps}{2}\!=\!1-\eps$, Algorithm \ref{algor:primes}
picks an $i$ with $M_i$ not divisible by any element of $\cE(x)$ and
a $c$ such that $p\!:=\!1+cM_i$ is prime. In
particular, we clearly have that 
$\log p \!=\!O(\log(1+KM_i))\!=\!O(n\log(n)+\log(s/\eps))$. \qed

\medskip
\noindent
{\bf (Complexity Analysis of Algorithm \ref{algor:primes}):} Let $L'\!:=\!nL$
and, for the remainder of our proof, let $p_{i}$ denote the $i\thth$ prime.
Since $L'\!\geq\!6$, $p_{L'}\!\leq L'(\log(L') + \log \log L')$ by
\cite[Thm.\ 8.8.4, pg.\ 233]{bs}. Recall that the
primes in $[\cL]$ can be listed simply by deleting all multiples of
$2$ in $[\cL]$, then deleting all multiples of $3$ in $[\cL]$,
and so on until one reaches multiples of $\lfloor \sqrt{\cL}\rfloor$.
(This is the classic sieve of Eratosthenes.) Recall also that one
can multiply an integer in $[\mu]$ and an integer $[\nu]$
within $O((\log\mu)(\log \log\nu)(\log\log\log \nu)
+(\log\nu)(\log \log\mu) \log\log\log \mu)$ bit operations
(see, e.g., \cite[Table 3.1, pg.\ 43]{bs}). So let us
define the function $\lambda(a):=(\log\log a)\log\log\log a$.

\noindent
{\bf Step 0:} By our preceding observations, it is easily checked that
Step 0 takes $O(L'^{3/2}\log^3 L')$ bit operations.

\noindent
{\bf Step 1:} This step consists of $n-1$ multiplications of
primes with $O(\log L')$ bits (resulting in $M_L$, which has
$O(n\log L')$ bits),
multiplication of a small power of $M_L$ by a square root of $M_L$,
division by an integer with $O(n\log L')$ bits, a constant number of 
additions of integers of
comparable size, and the generation of $O(\log L)$ random bits.
Employing Remark \ref{rem:tran} along the way, we thus arrive routinely at
an estimate of \\
\mbox{}\hfill
$O\left(n^2(\log L')\lambda(L')+\log(1/\eps)\lambda(1/\eps))
\right)$ \hfill \mbox{}\\
for the total number of bit operations needed for Step 1.

\noindent
{\bf Step 2:} Similar to our analysis of Step 1, we see that
Step 2 has bit complexity\\
\mbox{}\hfill $O((n\log(L')+\log(1/\eps))\lambda(n\log L'))$. \hfill
\mbox{}

\noindent
{\bf Step 3:} This is our most costly step: Here, we require\\
\mbox{}\hfill $O(\log K)\!=\!O(n\log(L')+\log(1/\eps))$\hfill\mbox{}\\
random bits and $J\!=\!O(\log x)\!=\!O(n\log(L')+\log(1/\eps))$ primality tests
on integers with
$O(\log(1+cM_i))\!=\!O(n\log(L')+\log(1/\eps))$ bits.
By an improved version of the AKS primality testing algorithm \cite{aks,lp}
(which takes $O(N^{6+\delta})$ bit operations to test an $N$ bit integer for
primality), Step 3 can then clearly be done within\\
\mbox{}\hfill
$O\!\left((n\log(L')+\log(1/\eps))^{7+\delta}\right)$
\hfill\mbox{}\\
bit operations, and the generation of $O(n\log(L')+\log(1/\eps))$ random
bits.

\noindent
{\bf Step 4:} This step clearly takes time on the order of the number of
output bits, which is just $O(n\log(n)+\log(1/\eps))$ as already observed
earlier.

\medskip
\noindent
{\bf Conclusion:} We thus see that Step 0 and
Step 3 dominate the complexity of our algorithm, and we are left with an
overall randomized complexity bound of\\
\mbox{}\hfill
$O\!\left(L'^{3/2}\log^3(L')+ \left(n\log(L')+\log(1/\eps)\right)^{7+\delta}
\right)$\hfill\mbox{}\\
\mbox{}\hfill $=O\!\left(\left(\frac{n}{\eps}\right)^{3/2}\log^3(n/\eps)
+\left(n\log(n)
+\log(1/\eps) \right)^{7+\delta}\right)$\hfill\mbox{}\\
\mbox{}\hfill $=O\!\left(\left(\frac{n}{\eps}\right)^{\frac{3}{2}+\delta}+
\left(n\log(n)+\log(1/\eps)\right)^{7+\delta}
\right)$\hfill\mbox{}\\
randomized bit operations. \qed

\subsection{Transferring from Complex Numbers to p-adics}  
\label{sub:transfer} 
\begin{prop}
\label{prop:sos} 
Given any $f_1,\ldots,f_k\!\in\!\Z[x_1]$ with  
maximum coefficient absolute value $H$, let $d_i\!:=\!\deg f_i$ and \\
\mbox{}\hfill $\tilde{f}(x_1):=
x^{d_1}_1f_1(x_1)f_1(1/x_1)+\cdots+x^{d_k}_1f_k(x_1)f_k(1/x_1)$.\hfill 
\mbox{}\\ 
Then $f_1\!=\cdots=\!f_k\!=\!0$ has a root on the complex unit circle iff 
$\tilde{f}$ has a root on the complex unit circle. In particular, if 
$f_i\!\in\!\cF_{1,\mu_i}$ and $\mu_i\!\leq\!m$ for all $i$, then 
$\tilde{f}\!\in\!\cF_{1,\mu}$ for some $\mu$ with $\mu\!\leq\!((m-1)m+1)k$ 
and $\tilde{f}$ has maximum coefficient bit-size $O(\log(kmH))$. \qed 
\end{prop} 

\noindent
Proposition \ref{prop:sos} follows easily upon observing that\linebreak  
$f_i(x_1)f_i(1/x_1)\!=\!|f_i(x_1)|^2$ for all $i\!\in\![k]$ and any $x_1 
\!\in\!\C$ with $|x_1|\!=\!1$. 

\begin{lemma}
\label{lemma:syl} 
\scalebox{.87}[1]{(See, e.g., \cite[Ch.\ 12, Sec.\ 1, pp.\ 397--402]{gkz94}.)}
\linebreak   
Suppose $f(x_1)\!=\!a_0+\cdots+a_dx^{d}_1$ and 
$g(x_1)\!=\!b_0+\cdots+b_{d'}x^{d'}_1$ are polynomials with 
indeterminate coefficients. Define their {\bf Sylvester matrix} 
to be the $(d+d')\times (d+d')$ matrix 

\noindent
\mbox{}\hfill \scalebox{.8}[.8]{$\cS_{(d,d')}(f,g)\!:=\!\begin{bmatrix}
a_0 & \cdots & a_d    & 0       & \cdots & 0 \\
   & \ddots &  & & \ddots &  \\
0   & \cdots & 0 & a_0 & \cdots & a_d \\
b_0 & \cdots & b_{d'}    & 0       & \cdots & 0 \\
  & \ddots &  &  & \ddots &  \\
0   & \cdots & 0 & b_0 & \cdots & b_{d'} 
\end{bmatrix}
\begin{matrix}
\\
\left. \rule{0cm}{.9cm}\right\}
d' \text{ rows}\\
\left. \rule{0cm}{.9cm}\right\}
d \text{ rows} \\
\\
\end{matrix}$}\hfill\mbox{}\\
\scalebox{.87}[1]{and their {\bf Sylvester resultant} to be 
$\cR_{(d,d')}(f,g)\!:=\!\det \cS_{(d,d')}(f,g)$.}\linebreak  
Then, assuming $f,g\!\in\!K[x_1]$ for some field $K$ and $a_db_{d'}\!\neq\!0$, 
we have that $f\!=\!g\!=\!0$ has a root in the algebraic closure of 
$K$ iff $\cR_{(d,d')}(f,g)\!=\!0$. 
Finally, if we assume further that $f$ and $g$ have complex coefficients of 
absolute value $\leq\!H$, and $f$ (resp.\ $g$) has exactly $m$ 
(resp.\ $m'$) monomial terms, then $|\cR_{(d,d')}(f,g)|\!\leq\!
m^{d'/2}m'^{d/2}H^{d+d'}$. \qed 
\end{lemma} 
\noindent
The last part of Lemma \ref{lemma:syl} follows easily from 
Hadamard's Inequality (see, e.g., \cite[Thm.\ 1, pg.\ 259]
{mignotte}). 

\begin{lemma}
\label{lemma:red} 
Suppose $D\!\in\!\N$ and $f\!\in\!\Z[x_1]\!\setminus\!\{0\}$ has degree $d$,   
exactly $m$ monomial terms, and maximum coefficient absolute value $H$. Also 
let $p$ be any prime congruent to $1$ mod $D$.  
Then $f$ vanishes at a complex $D\thth$ root of unity $\Longleftrightarrow 
f$ vanishes at a $D\thth$ root of unity in $\Q_p$. \qed 
\end{lemma} 
\begin{rem} 
Note that $x^2_1+x_1+1$ vanishes at a $3\rd$ root of unity in $\C$,  
but has {\bf no} roots at all in $\F_5$ or $\Q_5$. Hence our 
congruence assumption on $p$ in Lemma \ref{lemma:red}. \dia 
\end{rem} 

\noindent
{\bf Proof of Lemma \ref{lemma:red}:}
First note that by our assumption
on $p$, $\Q_p$ has $D$ distinct $D\thth$ roots of
unity: This follows easily from Hensel's Lemma (cf.\ the Appendix) 
and $\F_p$ having $D$ distinct $D\thth$ roots of unity.
Since $\Z\hookrightarrow\Q_p$ and $\Q_p$
contains all $D\thth$ roots of unity by construction, the
equivalence then follows directly from Lemma \ref{lemma:syl}.  \qed

\subsection{A Remark on Natural Density} 
Let us now introduce the {\bf $\pmb{\cA}$-discriminant} and clarify how 
often our $p$-adic speed-ups hold for inputs with bounded coefficients.
\begin{dfn}
Write any $f\!\in\!\C[x_1]$ as $f(x_1)\!=\!\sum^m_{i=1}c_ix^{a_i}_1$
with $0\!\leq\!a_1\!<\cdots<\!a_m$. Letting $\cA\!=\!\{a_1,\ldots,a_m\}$,
and\linebreak following the notation of Lemma \ref{lemma:red}, we
then define $\cD_\cA(f)$ to be 
\mbox{}\hfill $\cR_{(a_m-a_1,a_m-a_2)}\left.\left(\frac{f(x_1)}
{x^{a_1}_1},\left.\frac{\partial\left(\frac{f(x_1)}{x^{a_1}_1}\right)}
{\partial x_1}\right/x^{a_2-1}\right)\right/c_m$\hfill\mbox{}\\
to be the {\bf $\pmb{\cA}$-discriminant} of $f$ (see also
\cite[Ch.\ 12, pp.\ 403--408]{gkz94}). Finally, if $c_i\!\neq\!0$ for all $i$,
then we call $\supp(f)\!:=\!\{a_1,\ldots,a_m\}$ the {\bf support} of $f$. \dia
\end{dfn}
\begin{cor}
\label{cor:lots}
For any subset $\cA\!\subset\!\N\cup\{0\}$ of cardinality $m$,
let $\cI_\cA$ denote the family of pairs
$(f,p)\!\in\!\Z[x_1]\times \Pro$ with $f(x)\!=\!\sum^m_{i=1}c_ix^{a_i}_1$
and let $\cI^*_\cA$ denote the subset of $\cI_\cA$
consisting of those pairs $(f,p)$ with
$p\not\!|\cD_\cA(f)$. Also let $\cI_\cA(H)$ (resp.\ $\cI^*_\cA(H)$)
denote those pairs $(f,p)$ in $\cI_\cA$ (resp.\ $\cI^*_\cA$) where
$|c_i|\!\leq\!H$ for all $i\!\in\![m]$ and $p\!\leq\!H$. Then
$\frac{\#\cI^*_\cA(H)}{\#\cI_\cA(H)}\!\geq\!\left(1-\frac{(2d-1)m}{H}\right)
\left(1-\frac{d\log_2(dmH)}{H}\right)$. \qed
\end{cor}

\noindent
Our corollary above follows easily from our proof of
Assertion (3) of Theorem \ref{thm:qp} via an application of Lemma
\ref{lemma:syl} and the Schwartz-Zippel Lemma \cite{schwartz}, 
and is {\bf not} used in any of our proofs. 

\section{The Proof of Theorem 1.4}  
\label{sec:qp}   

\noindent          
{\bf (Assertion (1): $\pmb{\fqp(\cF_{1,m}\times\Pro)\!\in\!\pp}$ for 
$m\!\leq\!2$):} First note that the case $m\!\leq\!1$ is trivial: 
such a univariate $m$-nomial has no roots in $\Q_p$ iff  
it is a nonzero constant. So let us now assume $m\!=\!2$. 

Next, we can easily reduce 
to the special case $f(x)\!:=\!x^d-\alpha$ with $\alpha\!\in\!\Q$,  
since we can divide any input by a suitable monomial term, and arithmetic over 
$\Q$ is doable in polynomial time. The case $\alpha\!=\!0$ always results in 
the root $0$, so let us also assume $\alpha\!\neq\!0$. Clearly then,
any $p$-adic root $\zeta$ of $x^d-\alpha$ satisfies
$d\ord_p\zeta\!=\!\ord_p\alpha$. Since we can compute $\ord_p\alpha$
and reductions of integers mod $d$ in polynomial-time \cite[Ch.\ 5]{bs}, we 
can then assume that $d|\ord_p\alpha$ (for otherwise, $f$ 
would have no roots over $\Q_p$). Replacing $f(x_1)$ by 
$p^{-\ord_p\alpha}f(p^{\ord_p\alpha/d}x_1)$, we can assume further that 
$\ord_p\alpha\!=\!\ord_p\zeta\!=\!0$. In 
particular, if $\ord_p\alpha$ was initially a nonzero multiple of $d$, 
then $\log \alpha\!\geq\!d\log_2 p$. So $\size(f)\!\geq\!d$ and our rescaling 
at worst doubles $\size(f)$.  
                         
Letting $k\!:=\!\ord_p d$, note that $f'(x)\!=\!dx^{d-1}$ and thus 
$\ord_p f'(\zeta)\!=\!\ord_p(d)+(d-1)\ord_p\zeta\!=\!k$. So by Hensel's 
Lemma (cf.\ the Appendix), 
it suffices to decide whether
the $\mod \ p^\ell$ reduction of $f$ has a root in $(\Z/p^\ell\Z)^*$,
for $\ell\!=\!1+2k$. Note in particular that $\size(p^\ell)\!=\!
O(\log(p)\ord_p d)\!=\!O(\log(p)\log(d)/\log p)\!=\!O(\log d)$ which is 
linear in our notion of input size. By
Lemma \ref{lemma:qp} of the Appendix, we can then clearly decide whether 
$x^d-\alpha$ 
has a root in $(\Z/p^\ell\Z)^*$ within $\pp$ (via a single fast 
exponentiation), provided $p^\ell\!\not\in\!\{8,16,32,\ldots\}$.
                                                                                
To dispose of the remaining cases $p^\ell\!\in\!\{8,16,32,\ldots\}$, 
first note that we can replace $d$ by its reduction mod $2^{\ell-2}$ 
since every element of $(\Z/2^\ell\Z)^*$ has order dividing $2^{\ell-2}$,  
and this reduction can certainly be computed in polynomial-time. 
Let us then write $d\!=\!2^h d'$ where $2\!\!\not\!|d'$ and $h\!\in\!\{0,
\ldots,\ell-3\}$, and compute $d''\!:=\!1/d' 
\ \mod \ 2^{\ell-2}$. Clearly then, $x^d-\alpha$ has a root in 
$(\Z/2^\ell\Z)^*$ iff $x^{2^h}-\alpha'$ has a root in 
$(\Z/2^\ell\Z)^*$, where $\alpha'\!:=\alpha^{d''}$ (since 
exponentiation by any odd power is an automorphism of $(\Z/2^\ell\Z)^*$).  
Note also that $\alpha'$, $d'$, and $d''$ can clearly be computed 
in polynomial time. 

Since $x^{2^h}-\alpha'$ always has a root in $(\Z/2^\ell\Z)^*$  
when $h\!=\!0$, we can then restrict our root search to the cyclic 
subgroup $\{1,5^2,5^4,5^6,\ldots,5^{2^{\ell-2}-2}\}$ when $h\!\geq\!1$ 
and $\alpha'$ is a square (since there can be no roots when 
$h\!\geq\!1$ and $\alpha'$ is not a square). Furthermore, 
we see that $x^{2^h}-\alpha'$ can have no roots in 
$(\Z/2^\ell\Z)^*$ if $\ord_2\alpha'$ is odd. So, by rescaling $x$, 
we can assume further that $\ord_2\alpha'\!=\!0$, and thus that $\alpha'$ 
is odd. Now an odd $\alpha'$ is a square 
in $(\Z/2^\ell\Z)^*$ iff $\alpha'\!\equiv\!1 \; \mod \; 8$ 
\cite[Ex.\ 38, pg.\ 192]{bs}, and this can 
clearly be checked in $\pp$. So we can at last decide the existence of a root 
in $\Q_2$ for $x^d-\alpha$ in $\pp$: Simply combine fast exponentiation with 
Assertion 3 
of Lemma \ref{lemma:qp} again, applied to $x^{2^h}-\alpha'$ over the cyclic 
group $\{1,5^2,5^4,5^6,\ldots,5^{2^{\ell-2}-2}\}$. 

\medskip 
\noindent
{\bf (Assertion (2): $\pmb{\fqp(\cF_{1,3}\times\Pro)\!\in\!\pp}$ for non-flat  
$\pmb{\newt_p(f)}$):}  First note that
$x\!\in\!\Q_p\setminus\Z_p \Longleftrightarrow \frac{1}{x}\!\in\!p\Z_p$.
Letting $f^*(x)\!:=\!x^{\deg f}f(1/x)$ denote the reciprocal polynomial
of $f$, note that the set of $p$-adic rational roots of $f$ is simply the
union of the $p$-adic integer roots of $f$ and the reciprocals of the
$p$-adic integer roots of $f^*$. So we need only show we can detect 
roots in $\Z_p$ in $\pp$. 

As stated, Assertion (2) then follows directly from Lemma \ref{lemma:ai}. 

So let us now concentrate on extending polynomiality to some of our 
exceptional inputs: Writing $f(x)\!=\!c_1+c_2x^{a_2}+c_3x^{a_3}$ as before, 
let us consider the special case where 
$f\!\in\!\cF_{1,3}$ has a degenerate root in $\C_p$ and 
$\gcd(a_2,a_3)\!=\!1$. Note that we now allow $p$ to divide any 
number from\linebreak $\{a_2,a_3,a_3-a_2\}$. (It is 
easily checked that the collinearity condition fails for such 
polynomials since their $p$-adic Newton polygons are line 
segments.) The $\{0,a_2,a_3\}$-discriminant of $f$ then 
turns out to be $\Delta:=(a_3-a_2)^{a_3-a_2}a^{a_2}_2c^{a_3}_2-
(-a_3)^{a_3}c^{a_3-a_2}_1c^{a_2}_3$ 
(see, e.g., \cite[Prop.\ 1.8, pg.\ 274]{gkz94}). In particular, 
while one can certainly evaluate $\Delta$ with a small number
of arithmetic operations, the bit-size of $\Delta$ can be
quite large. However, we can nevertheless efficiently decide whether $\Delta$ 
vanishes for integer $c_i$ via {\bf gcd-free bases} 
(see, e.g., \cite[Sec.\ 2.4]{brs}). Thus, we can at least 
check whether $f$ has a degenerate root in $\C_p$ in $\pp$. 

Given an $f$ as specified, it is then easily checked that if $\zeta\!\in\!\C_p$ 
is a degenerate root of $f$ then 
the vector $[c_1,c_2\zeta^{a_2},c_3\zeta^{a_3}]$ must be a right null 
vector for the matrix $\begin{bmatrix}1 & 1 & 1 \\ 0 & a_2 & a_3
\end{bmatrix}$. In other words, $[c_1,c_2\zeta^{a_2},c_3\zeta^{a_3}]$ is a 
mutiple of $[\alpha, \beta,\gamma]$ for some integers $\alpha,\beta,\gamma$ 
with size polynomial in $\size(f)$. Via the extended Euclidean algorithm 
\cite[Sec.\ 4.3]{bs}, we can 
find $A$ and $B$ (also of size 
polynomial in $\size(f)$) with $Aa_2+Ba_3\!=\!1$. So then we obtain that\\ 
\mbox{}\hfill $\left(\frac{c_2\zeta^{a_2}}{c_1}\right)^A
\left(\frac{c_3\zeta^{a_3}}
{c_1}\right)^B\!=\!\frac{c^A_2c^B_3}{c^{A+B}_1}\zeta\!=\!
\left(\frac{\beta}{\alpha}\right)^A\left(\frac{\gamma}
{\alpha}\right)^B$.\hfill\mbox{}\\ 
In other words, $f$ has a rational root, and thus 
this particular class of $f$ always has $p$-adic rational roots. \qed

\medskip
\noindent
{\bf (Assertion (3): $\pmb{\fqp(\Z[x_1]\times \Pro)\!\in\!\np}$ for 
most\linebreak  
inputs):} 
Just as in our reduction from $\Q_p$ to $\Z_p$ in the beginning of our last 
proof, it is enough to show that, for most $f$, roots in $\Z_p$ admit 
succinct certificates. We can also clearly assume that $f$ is not 
divisible by $x_1$. 

Observe now that the $p$-adic valuations of all the roots of $f$ in $\C_p$
can be computed in polynomial-time. This is
easily seen via two facts: (1) convex hulls of subsets of $\Z^2$
can be computed in polynomial-time (see, e.g., \cite{edelsbrunner}), and
(2) the valuation of any
root of $f(x)\!=\!\sum^m_{i=1}c_ix^{a_i}$ is a ratio of the
form $\frac{\ord_p(c_i)-\ord_p(c_j)}{a_j-a_i}$, where
$(a_i,\ord_p(c_i))$ and $(a_j,\ord_p(c_j))$ are respectively the
left and right vertices of a lower edge of $\newt_p(f)$ (cf.\ Lemma
\ref{lemma:newt} of the Appendix).
Since $\ord_p(c_i)\!\leq\!\log_p(c_i)\!\leq\!\size(c_i)$, note in particular
that every root $\zeta\!\in\!\C_p$ of $f$
satisfies $|\ord_p\zeta|\!\leq\!2\max_i\size(c_i)\!\leq\!2\size(f)\!<\!
2\size_p(f)$.

Since $\ord_p(\Z_p)\!=\!\N\cup\{0\}$, we can clearly assume that
$\newt_p(f)$ has an edge with non-positive integral slope, for otherwise $f$
would have no roots in $\Z_p$. Letting $a$ denote the smallest nonzero
exponent in $f$, $g(x)\!:=\!f'(x)/x^{a-1}$, and $\zeta\!\in\!\Z_p$
any $p$-adic integer root of $f$, note then that
$\ord_p f'(\zeta)\!=\!(a-1)\ord_p(\zeta)+\ord_p g(\zeta)$. Note also that\\
\mbox{}\hfill $\cD_\cA(f)\!=\!\res_{a_m,a_m-a_1}(f,g)$ \hfill \mbox{}\\
so if $p\not\!\!|\cD_\cA(f)$ then
$f$ and $g$ have no common roots in the algebraic closure of
$\F_p$ by Lemma \ref{lemma:syl}. In particular,
$p\!\!\not|\cD_\cA(f)\Longrightarrow g(\zeta)\!\not
\equiv\!0 \; \mod \; p$; and thus $p\!\!\not\!|\cD_\cA(f,g)\Longrightarrow
\ord_p f'(\zeta)\!=\!(a-1)\ord_p(\zeta)$. Furthermore, by the convexity
of the lower hull of $\newt_p(f)$, it is clear that $\ord_p(\zeta)\!\leq\!
\frac{\ord_p c_i -\ord_p c_0}{a_1}\!\leq\!
\frac{2\max_i \log_p|c_i|}{a_1}$. So $p\not\!|\cD_\cA(f)\Longrightarrow
\ord_p f'(\zeta)\!<\!2\size(f)$. 

Our fraction of inputs admitting a succinct certificate will then
correspond precisely to those $(f,p)$ such that $p\!\!\not\!|\cD_\cA(f)$.
In particular, let us define $E$ to be the union of all pairs
$(f,p)$ such that $p|\cD_\cA(f)$, as $\cA$ ranges over all finite
subsets of $\N\cup\{0\}$. It is then easily checked that $E$ is a
countable union of hypersurfaces.

Fix $\ell\!=\!4\size(f)$. Clearly then,
by Hensel's Lemma, for any $(f,p)\!\in\!(\Z[x_1]\times \Pro)\setminus E$,
$f$ has a root $\zeta\!\in\!\Z_p \Longleftrightarrow f$ has a root
$\zeta_0\!\in\!\Z/p^\ell\Z$. Since
$\log(p^\ell)\!=\!O(\size(f)\log p)\!=\!O(\size_p(f)^2)$,
and since arithmetic in $\Z/p^\ell\Z$
can be done in time polynomial in $\log(p^\ell)$ \cite[Ch.\ 5]{bs}, we have
thus at last found our desired certificate: a root
$\zeta_0\!\in\!(\Z/p^\ell\Z)^*$ of $f$ with $\ell\!=\!4\size(f)$.

To conclude, the assertion on checking whether trinomial inputs 
lie in $E$ follows immediately from our earlier observations on deciding 
the vanishing of $\Delta$. In particular, instead of applying gcd-free 
bases, we can instead simply use recursive squaring and efficient 
$\F_p$-arithmetic. \qed 

\medskip 
\noindent 
{\bf (Assertion (4): $\pmb{\fqp(\Z[x_1]\times \Pro)}$ is $\np$-hard \linebreak 
under $\zpp$-reductions):}  
We will prove a ($\zpp$) randomized polynomial-time reduction from $\sat$ to
\linebreak 
$\fqp(\Z[x_1]\times\Pro)$, making use of the intermediate input families
$\{(\Z[x_1])^k\; | \; k\!\in\!\N\}$ and $\Z[x_1]\times\{x^D_1-1\;
| \; D\!\in\!\N\}$ along the way.

Toward this end, suppose $B(y)\!:=\!C_1(y)\wedge\cdots\wedge C_k(y)$
is any $\sat$ instance. The polynomial system
$(\cP_P(C_1),\ldots,$\linebreak
$\cP_P(C_k))$, for
$P$ the first $n$ primes (employing Lemma \ref{lemma:plai}), then
clearly yields the implication\linebreak $\feas_\C(\{(\Z[x_1])^k\; | \;
k\!\in\!\N\})\!\in\!\pp \Longrightarrow \pp\!=\!\np$.
Composing this reduction with Proposition \ref{prop:sos}, we then
immediately obtain the implication $\feas_\C(\Z[x_1]\times\{x^D_1-1\; | \;
D\!\in\!\N\})\!\in\!\pp\Longrightarrow \pp\!=\!\np$.

At this point, we need only find a means of transferring from $\C$ to
$\Q_p$. This we do by preceding our reductions above by a judicious 
(possibly new) 
choice of $P$. In particular, by applying Theorem \ref{thm:von} with
$\eps\!=\!1/3$ (cf.\ Lemma \ref{lemma:red}) we 
immediately obtain the implication\linebreak 
$\fqp((\Z[x_1]\times\{x^D_1-1\; | \;
D\!\in\!\N\})\times \Pro)\!\in\!\zpp\Longrightarrow \np\!\subseteq\!\zpp$.

To conclude, observe that any root $(x,y)\!\in\!\Q^2_p\setminus\{(0,0)\}$ 
of the quadratic form $x^2-py^2$ must satisfy \linebreak 
$2\ord_p x\!=\!1+2\ord_p y$ 
--- an impossibility. Thus the only $p$-adic rational root of $x^2-p y^2$ is 
$(0,0)$ and we easily obtain a polynomial-time reduction from\linebreak  
$\fqp((\Z[x_1]\times\{x^D_1-1\; | \;
D\!\in\!\N\})\times\Pro)$ to\linebreak 
$\fqp(\Z[x_1]\times\Pro)$: simply map any 
instance\linebreak  
$(f(x_1),x^D_1-1,p)$ of the former problem to\linebreak  
$(f(x_1)^2-(x^D_1-1)^2p,p)$. So we are done. \qed

\medskip 
\noindent 
{\bf (Assertion (5): $\pmb{\fqp(\Z[x_1]\times \Pro)}$ is $\np$-hard,\linebreak 
assuming Wagstaff's Conjecture):}  
If we also have the truth of the Wagstaff Conjecture 
then we simply repeat our last proof, replacing our 
AGP Theorem-based algorithm with a simple brute-force search. 
This maintains polynomial complexity, but with the added advantage of 
completely avoiding randomization. \qed 

\section*{Acknowledgements} 
The authors would like to thank David Alan Plaisted for his kind 
encouragement, and Eric Bach, Sidney W.\ Graham, and Igor Shparlinski 
for many helpful comments on primes in arithmetic progression. We also 
thank Matt Papanikolas for valuable $p$-adic discussions. Finally, 
we thank an anonymous referee for insightful comments that greatly 
helped clarify our presentation. 

\bibliographystyle{acm}

\section{Appendix: Additional Background} 
Let us first recall briefly the following complexity classes 
(see also \cite{papa} for an excellent textbook treatment): 
\begin{itemize}
\item[$\pp$]{ The family of decision problems which can be done within time
polynomial in the input size.\footnote{Note that the underlying polynomial
depends only on the problem in question (e.g., matrix inversion, shortest path
finding, primality detection) and not the particular instance of the
problem.}}
\item[$\zpp$]{ The family of decision problems admitting a
randomized polynomial-time algorithm giving a correct answer,  
or a report of failure, the latter occuring with probability 
$\leq\!\frac{1}{2}$. }
\item[$\np$]{ The family of decision problems where a ``{\tt Yes}'' answer can
be {\bf certified} within time polynomial in the input size.}
\item[$\expt$]{ The family of decision problems solvable 
within time exponential in the input size.}
\end{itemize}

The classical Hensel's Lemma can be phrased as follows. 
\begin{lemma}
\label{lemma:hensel} 
\cite[Pg.\ 48]{robert} 
Suppose $f\!\in\!\Z_p[x_1]$ and $\zeta_0\!\in\!\Z_p$
satisfies $f(\zeta_0)\!\equiv\!0 \ (\mod \ p^\ell)$ and
$\ord_p f'(\zeta_0)\!<\!\frac{\ell}{2}$. Then there is
a root $\zeta\!\in\!\Z_p$ of $f$ with $\zeta\!\equiv\!\zeta_0 \
(\mod \ p^{\ell-\ord_p f'(\zeta_0)})$ and $\ord_p f'(\zeta)\!=\!\ord_p
f'(\zeta_0)$. \qed
\end{lemma}

The final tool we will need is a standard lemma on binomial equations
over certain finite groups. Recall that for any ring $R$, we denote its unit
group by $R^*$.
\begin{lemma}
\label{lemma:qp}
(See, e.g., \cite[Thm.\ 5.7.2 \& Thm.\ 5.6.2, pg.\ 109]{bs})
Given any cyclic group $G$, $a\!\in\!G$, and an integer $d$, the following
3 conditions are equivalent:\\
\mbox{}\hspace{1cm}1.\ the equation $x^d\!=\!a$ has a solution $a\!\in\!G$.\\
\mbox{}\hspace{1cm}2.\ the order of $a$ divides $\frac{\#G}{\gcd(d,\#G)}$.\\
\mbox{}\hspace{1cm}3.\ $a^{\#G/\gcd(d,\#G)}\!=\!1$.\\
Also, $\F^*_q$ is cyclic for
any prime power $q$, and $(\Z/p^\ell\Z)^*$ is cyclic
for any $(p,\ell)$ with $p$ an odd prime or $\ell\!\leq\!2$.
Finally, for $\ell\!\geq\!3$, $(\Z/2^\ell\Z)^*\!=\!
\{\pm 1,\pm 5,\pm 5^2,\pm 5^3,\ldots,\pm 5^{2^{\ell-2}-1} \ \mod \ 2^\ell\}$. \qed
\end{lemma}

\end{document}